\documentclass{jloganal}

\usepackage{amsmath,amssymb,bbm,amsthm}
\usepackage{comment}

\newcommand{\assign}{:=}

\newcommand{\um}{-}

\newcommand{\reals}{\mathbbm{R}}
\newtheorem{theorem}{Theorem}[section]

\newtheorem{proposition}[theorem]{Proposition}
\newtheorem{corollary}[theorem]{Corollary}

\begin{document}

\title{A constructive proof of Simpson's Rule}
\author{Thierry Coquand}
\address{Thierry Coquand\\
Computing Science Department, G\"oteborg University}
\author{Bas Spitters}
\address{Bas Spitters\\
Institute for Computing and Information Science, Radboud University Nijmegen}
\received{}
\revised{}

\begin{abstract}
 For most purposes, one can replace the use of Rolle's theorem and the mean
value theorem, which are not constructively valid, by the law of bounded
change~\cite{BridgerStolzenberg}.
The proof of two basic results in numerical analysis,
the error term for Lagrange interpolation and Simpson's rule, however seem to
require the full strength of the classical Rolle's Theorem. The goal of this
note is to justify these two results constructively, using ideas going back to
Amp\`ere~\cite{Ampere} and Genocchi~\cite{Genocchi}.
\end{abstract}
\maketitle

\section{Introduction}
Rolle's Theorem states that
if $f$ is differentiable on a real interval $[a,b]$ with $a<b$ and $f(a) =
f(b)=0$ then there exists
$c$ in $]a,b[$ such that $f'(c) = 0$. It implies directly the mean value theorem
that
if $f$ is differentiable on a real interval $[a,b]$ with $a<b$, then there
exists
$c$ in $]a,b[$ such that $f(b) - f(a) = (b-a)f'(c)$.
This is a key result in most text books in Analysis. It is rather direct to see
that it does not hold
constructively~\cite{Bishop/Bridges:1985}, and is replaced in this context by an
approximated form: if $f$ is differentiable on a real interval $[a,b]$ with
$a<b$ then for any $\epsilon>0$ there exists
$c$ in $]a,b[$ such that $|f(b) - f(a) - (b-a)f'(c)|<\epsilon$.  This more
complex formulation can be thought to be a problem of constructive mathematics.

 It can be argued however~\cite{BridgerStolzenberg} that most applications of
the mean value theorem can be replaced by the law of bounded change that we have
$|f(b)-f(a)|\leqslant M(b-a)$ if $f$ is uniformly derivable on $[a,b]$
and $|f'(x)|\leqslant M$ for all $x$ in $[a,b]$. The law of bounded change is
constructively valid and
a presentation based on the law of bounded change~\cite{BridgerStolzenberg}
appears as elegant as the classical treatment of basic analysis results.
Interestingly, a criticism of the mean value theorem, which could have been
written by a constructive mathematician, appears in Dieudonn\'e {\em Foundations
of Modern Analysis}~\cite{Dieudonne}:
``the trouble with that classical
formulation is that $\dots$ it completely conceals the fact that {\em nothing}
is known on the number $c$, except
that it lies between $a$ and $b$, and for most purposes, all one needs to know
is that $f'(c)$ is a number
which lies between that g.l.b. and l.u.b. of $f'$ in the interval $[a,b]$ (and
\emph{not} the fact that it actually is a value of $f'$).'' The goal of this
note is to analyse two results that at first seem to require the full
strength of the classical mean value theorem. The second of these results is
Simpson's rule for approximating an integral, which is indeed proved in the
exercises of~\cite{Dieudonne} using Rolle's Theorem. We show, using some ideas
going back to Amp\`ere~\cite{Ampere} and Genocchi~\cite{Genocchi}, how to
justify these results constructively.

 This note is organised as follows. We first present the two results we want to
analyse: Lagrange error formula and Simpson's rule. We then explain Genocchi's
formula, in a way which stresses the connection with the
work of Bridger and Stolzenberg~\cite{BridgerStolzenberg},
and show how it can be used instead of Rolle's Theorem.


\section{Lagrange error formula and Simpson's rule}

 In this section, we present two basic results, Theorems~\ref{lagrange} and
\ref{simpson}, and explain how they can be classically derived using the
following generalization of Rolle's Theorem.

\begin{theorem}[Classical generalized Rolle's theorem]\label{GenRolle} Let $f$
be $n$ times differentiable and have $n + 1$ zeroes in an interval $[a, b]$.
  Then $f^{(n)}$ has a zero in $[a, b]$.
\end{theorem}

\begin{proof}
If $f$ has $n+1$ zeroes $x_0<x_1<\dots < x_n$ then using Rolle's Theorem
$f'$ will have $n$ zeroes $y_0,\dots, y_{n-1}$ with $y_i$ in $]x_i,x_{i+1}[$.
Hence we obtain the result by induction on $n$.
\end{proof}



We now present the Lagrange polynomial as it can be found in numerical
analysis textbooks; e.g.~\cite{BurdenFaires}.

\begin{theorem}[Lagrange error formula]\label{lagrange}
Let $f$ be $n$ times differentable on an interval $[a,b]$, $P$ the polynomial of
degree
$n-1$ which agrees with $f$ on $n$ values $x_0<\dots<x_{n-1}$ and $M$ such that
$|f^{(n)}(x)|\leqslant M$ for
all $x$ in $[a,b]$. Then for all $x$ in $[a,b]$
  \[ |f (x) - P (x) | \leqslant \frac{|\prod (x - x_k)|}{n!}M
  \]
\end{theorem}

\begin{proof}
 This is proved using Rolle's Theorem in the following
way~\cite{Dieudonne,BurdenFaires}.
First, classically, we can assume that $x$ is not equal to one of the $x_i$
since the inequality is
clear if $x=x_i$. We then consider the function
$$
g(y) = \frac{\prod (x - x_k)}{n!}(f(y)-P(y)) - \frac{\prod (y -
x_k)}{n!}(f(x)-P(x))
$$
This function is $n$ times differentiable and has $n+1$ zeroes
$x,x_0,\dots,x_{n-1}$.
Using Theorem \ref{GenRolle}, there exists $c$ such that $g^{(n)}(c) = 0$ which
can be written as
$$
f(x) - P(x) = \frac{\prod (x - x_k)}{n!}f^{(n)}(c).
$$
This finishes the proof.
\end{proof}



 By a similar use of Theorem~\ref{GenRolle}, one can derive the following
classical result~\cite{Dieudonne,BurdenFaires}.

\begin{theorem}[Simpson's rule~\cite{BurdenFaires}]\label{simpson}
If $f$ is $4$-differentiable on an interval $[a,b]$ and
$|f^{(4)}(x)|\leqslant M$ for all $x$ in $[a,b]$, then we have
  \[ \left| \int_a^b f (x) \hspace{0.25em} d x - \frac{b - a}{6} \left[ f (a)
     + 4 f \left( \frac{a + b}{2} \right) + f (b) \right] \right| \leqslant
     \frac{(b - a)^5}{2880} M. \]

\end{theorem}

Rolle's theorem \emph{is} constructively provable~\cite{Bishop/Bridges:1985}
provided $f'$ is
\emph{locally non-constant}: In every interval there are $x,y$ such that
$f'(x)<f'(y)$. In this case, $f'$ is in particular locally nonzero, i.e.\ in every interval there
is an $x$ such that $f'(x)>0$ or $f'(x)<0$. It follows that $f$ is locally non-constant, as
is readily seen by integration. We see that if $f^{(n)}$ is locally nonzero,
then for all $k<n$, $f^{(k)}$ is locally nonconstant. We obtain an equal
conclusion version of the generalized Rolle's theorem: Let $f$ be $n$ times
differentiable and have $n + 1$ zeroes in an interval $[a, b]$. If, moreover,
$f^{(n)}$ is locally nonzero, then $f^{(n)}$ has a zero in $[a, b]$.

From this equal conclusion version, we can obtain an equal hypothesis version
of Rolle's theorem.

\begin{proposition}
  Let $f$ be $n$ times differentiable and have $n + 1$ zeroes, $x_i$, in an
interval $[a, b]$. Then there exists $x \in [a, b]$ such that $|f^{(n)} (x)|
\leq \varepsilon$.
\end{proposition}

\begin{proof}
Either $\inf_{x \in [a, b]}|f^{(n)} (x) |<\varepsilon$ or $\inf_{x \in [a, b]}|f^{(n)} (x) |
> \varepsilon/2$. In the former case, we are done. In the latter case,
$f^{(n)}$ is locally non-zero, hence by the remark above, we can follow the
proof of Theorem~\ref{GenRolle} to conclude that $f^{(n)}$ has a zero. A
contradiction.
\end{proof}

In this way we can derive the Lagrange error formula~\ref{lagrange}. One can argue
however that this derivation is more complex than the classical result.
In Section~\ref{lagrange-pf} we give smooth proofs of Theorems~\ref{lagrange}
and~\ref{simpson}, which hold both classically {\em and} constructively,
and do not rely on Theorem~\ref{GenRolle}.

\section{\label{sec:HG}Hermite-Genocchi formula}
 The definition in~\cite{BridgerStolzenberg} of uniform differentiability of a
function $f$ on an interval $I = [a,b]$ can be formulated as follows: there
exists an uniformly continuous
function $F:I^2\rightarrow\reals$
such that $f(y) - f(x) = (y-x)F(x,y)$ for all $x,y$ in $I$. We can then define
$f'(x)$ to be $F(x,x)$.
The following result is a generalization of this characterisation
to $n$-differentiability. The proof is a simple application
of the fundamental theorem of the calculus~\cite{BridgerStolzenberg}.

\begin{theorem}\label{main}
A function $f:I\rightarrow \reals$ is uniformly $n$-differentiable if, and only
if, there exist $n+1$ uniformly
continuous functions
$f_0(x),~f_1(x_0,x_1),~\dots,~f_{n}(x_0,\dots,x_{n})$ defined respectively on
$I,~I^2,~\dots,~I^{n+1}$ and such that
\begin{eqnarray}
f_0(x)= f(x),~f_0(x_1)-f_0(x_0) = (x_1-x_0)f_1(x_0,x_1),\dots,\label{divdiff}\\
f_{n-1}(x_0,\dots,x_{n-2},x_{n})-f_{n-1}(x_0,\dots,x_{n-2},x_{n-1}) =
(x_{n}-x_{n-1})f_{n}(x_0,\dots,x_{n})\nonumber
\end{eqnarray}
for all $x_0,\dots,x_{n}$ in $I$.

We then have
$$
f^{(n)}(x) = n!{f_{n}(x,\dots,x)}
$$
and, conversely, we can define
\begin{equation}
f_n(x_0,\dots,x_n) = \int _{\Sigma_n} f^{(n)}(t_0x_0+\dots + t_nx_n) dt_0\dots
dt_n\label{Genocchi}
\end{equation}
where\footnote{%
This integral of a uniformly continuous function $g$ over $\Sigma_n$ can be
defined by induction on $n$ in the following way: for $n = 0$ the integral is
$g(1)$, and for $n>0$ we define
\[
 h(t_0,...,t_{n-1}) =  \int_0^1 g(t_0,...,t_{n-1}(1-u), t_{n-1}u) du
\]
and the integral of $g$ over $\Sigma_n$ is the integral of $h$ over
$\Sigma_{n-1}$.}
 $\Sigma_n = \{(t_0,\dots,t_n)\in [0,1]^{n+1}~|~t_0+\dots + t_n = 1\}$.
\end{theorem}
\begin{proof}
Assume that the divided differences are uniformly continuous on
$I,I^2,\dots,I^{n+1}$.
Since
$$f(y)-f(x) = (y-x)f_1(x,y)$$
we get that $f$ is uniformly differentiable and $f'(x) = f_1(x,x)$. We then
have
  \begin{eqnarray}
f'(y)-f'(x)  & = & f_1(y,y)-f_1(x,x)   \nonumber\\
& = & f_1(y,y)-f_1(y,x)+f_1(y,x)-f_1(x,x) \nonumber\\
& = &  (y-x)(f_2(y,x,y)+f_2(x,x,y)) \nonumber
  \end{eqnarray}
and so $f$ is uniformly $2$-differentiable and $f^{(2)}(x) = 2f_2(x,x,x)$.
Proceeding in this way
we see that $f$ is uniformly $k$-differentiable and that we have $k!
f_k(y,\dots,y) = {f^{(k)}(y)}$ for
$k=0,\dots,n$.

 Conversely, assume that $f$ is uniformly $n$-differentiable. We define
$$
f_k(x_0,\dots,x_k) = \int _{\Sigma_k} f^{(k)}(t_0x_0+\dots + t_kx_k) dt_0\dots
dt_k
$$
for $k = 0,\dots,n$. These functions are uniformly continuous. Furthermore, we
have $f_0(x) = f(x)$ and
$$
f_0(x)-f_0(x_0) = (x-x_0)f_1(x_0,x)
$$
since
$$
(x-x_0)f_1(x_0,x) = \int_0^1 (f((1-t)x_0+tx))' dt
$$
holds for $x$ apart from $x_0$, by the fundamental theorem of the
calculus~\cite{BridgerStolzenberg} and hence
for all $x,x_0$ by continuity. It follows that we have
  \begin{eqnarray}
& & f_k(x_0,\dots,x_{k-1},x) - f_k(x_0,\dots,x_{k-1},x_k)    \nonumber\\
& = &
\int_{\Sigma_{k}} (f^{(k)}(t_0x_0+\dots + t_{k}x) - f^{(k)}(t_0x_0+\dots
+ t_{k}x_k)) dt_0\dots dt_{k}
    \nonumber\\
& = & (x-x_k) \int_{\Sigma_{k}} \int_0^1 (f^{(k+1)}(t_0x_0+\dots + t_{k}(1-u)x_k+t_k ux) du dt_0\dots dt_{k}
    \nonumber\\
& = & (x-x_k)\int_{\Sigma_{k+1}} f^{(k+1)}(t_0x_0+\dots + t_kx_k + t_{k+1}x) dt_0\dots dt_{k+1}
    \nonumber\\
& = & (x-x_k) f_{k+1}(x_0,\dots,x_k,x) \nonumber
  \end{eqnarray}
which shows that these functions satisfy the required equations.
\end{proof}


Following Ampere~\cite{Ampere} we observe from Formula~\ref{divdiff}:
\begin{eqnarray}
 f_1(x_0,x_1) &=&   \frac{f(x_1)}{x_1-x_0}  + \frac{f(x_0)}{x_0-x_1}\nonumber\\
 f_2(x_0,x_1,x_2) &=& \frac{f(x_2)}{(x_2-x_0)(x_2-x_1)} + \frac{f(x_1)}{(x_1-x_0)(x_1-x_2)} +
\frac{f(x_0)}{(x_0-x_1)(x_0-x_2)}\nonumber\\
f_3(x_0,x_1,x_2,x_3) &= & \ldots\nonumber
\end{eqnarray}
when all he $x_i$'s are distinct. Hence, the functions $f_k(x_0,\dots,x_k)$ are symmetric --- that is, they are
invariant under permutation of the variables. Also by formula~\ref{divdiff},
$$
f(x) = f_0(x_0) + (x-x_0)f_1(x_0,x_1)+ \dots +
        (x-x_0)\dots (x-x_{n-1})f_{n}(x_0,\dots,x_{n-1},x).
$$
Formula~\ref{Genocchi} is known as Hermite-Genocchi formula. Genocchi found these formulae by analysing
the notion of ``fonctions interpolaires''
due to Amp\`ere~\cite{Ampere}\footnote{In the case where $f$ is a monic
polynomial,
the functions $f_1(x_0,x_1),~\dots,~f_n(x_0,\dots,x_{n-1})$ are also polynomial
and they form with $f$ a Gr\"obner basis of the universal decomposition
algebra of the polynomial $f$; see e.g.~\cite{rennert1999calcul}.}.


 A direct consequence of Theorem \ref{main} is the following result.

\begin{corollary}
Given a function $f:I\rightarrow \reals$ and $n+1$ distinct elements
$x_0,\dots,x_n$ in $I$,		 let $P(x_0,\dots,x_{n},x)$ be the interpolation
polynomial
$a_{n-1}(x_0,\dots,x_n)x^{n-1} +\dots+a_0(x_0,\dots,x_n)$
of $f$ at $x_0,\dots,x_n$. Then $P(x_0,\dots,x_n,x)$, seen
as a function of the parameters $x_0,\dots,x_n$, can be extended to an uniformly
continuous function on $I^{n+1}$ (that is, each function $a_i(x_0,\dots,x_n)$
can be extended to an uniformly continuous function on $I^{n+1}$) if, and only
if,
$f$ is uniformly $n$-differentiable.
\end{corollary}


\section{Applications}

 We explain now how to derive Theorem \ref{lagrange} from Theorem \ref{main}.
We assume that $f:I\rightarrow\reals$ is uniformly $n$-differentiable.
Given any $n$ elements $x_0,\dots,x_{n-1}$ in $I$ we associate the
\emph{Newton polynomial} of degree $n-1$:
\[ P (x) : = f_0(x_0) + f_1(x_0, x_1) (x - x_0) + \cdots + f_{n-1}(x_0, ..., x_{n-1}) (x
   - x_0) (x - x_1) \cdots (x - x_{n - 2}) . \label{Newton-coeff} \]
We have
\[ f(x) - P(x) = (x-x_0)\dots (x-x_{n-1})f_n(x_0,\dots,x_{n-1},x) \]
On the other hand, we also have $P(x_i) = f(x_i)$ for $i=0,\dots,n-1$.

If $|f^{(n)}(u)|\leqslant M$ for all $u$ in $I$, then by Theorem~\ref{main},
\[ |f_n(x_0,\dots,x_{n-1},x)|\leqslant
M\int _{\Sigma_n} 1 = \frac{M}{n!} \]
The last equality follows from applying Formula~\ref{Genocchi} to the function $g(x):=\frac{x^n}{n!}$ and
observing that $g^{(n)}=1$ and $g_n(x,\ldots,x)=\frac1{n!}$, by Formula~\ref{divdiff}. This proves
Theorem~\ref{lagrange}.

 Notice that the Newton polynomial $P(x)$ is defined for any $x_0,\dots,x_{n-1}$
without requiring them to be distinct. We have $P^{(l-1)}(x_i) = f^{(l-1)}(x_i)$
if $x_i$ is duplicated $l$ times in this list.

Theorem~\ref{main} is also valid for $f : I \to E$, where $E$ is a  Banach space.
Theorem~\ref{lagrange} is also valid for $f : I \to F$, where $F$ is a normed space.
See for instance~\cite{Dieudonne,Banach-algs} for differentation and integration with values in a Banach space.

\subsection{A constructive proof of Simpson's
rule}\label{lagrange-pf}
As explained above, the typical proofs of Simpson's rule~\ref{simpson},
see e.g.~\cite{BurdenFaires, SuliMayers}, use Rolle's theorem, and so are not
constructively valid.
We adapt the proof in~{\cite{SuliMayers}} which uses Rolle's theorem
three times on the triple zero at 0 and the simple zero at 1 of the function $H$
below.

Define $F (t) \assign f ( \frac{a + b}{2} + \frac{b - a}{2} t)$. This reduces
the problem to showing that
\[| \int_{- 1}^1 F (\tau) d \tau - \frac{1}{3} (F (- 1) + 4 F (0) + F (1))
| \leqslant N / 90,\]
where $N \assign \|F^{(4)} \|.$

  Define
\[ G (t) = \int_{- t}^t F (\tau) d\tau -\frac{t}{3} (F(-t) + 4 F(0) + F(t))\]
  We need to prove that $90 G (1) \leqslant N$. To do so, define $H (t)
  \assign G (t) - t^5 G (1)$. Then
  \[ H (0) = H (1) = H' (0) = H'' (0) = 0. \]
  Hence, $H_3(0, 0, 0, 1) = - (H_2(0, 0, 0) - H_2(0, 0, 1)) = 0 + (\um H _1(0, 0) +
  H_1(0, 1)) = 0$. (This line replaces three uses of Rolle's theorem.)

  Moreover,\[H^{(3)} (t) = - \frac{t}{3} (F^{(3)} (t) - F^{(3)} (- t)) - 60 t^2 G (1) = - \frac{t}{3} ( \int_{- t}^t F^{(4)}) - 60 t^2 G (1).\]

  This shows that
  \begin{eqnarray}
    0 = H (0, 0, 0, 1) & = & \int^1_0 H^{(3)} \nonumber\\
    & = & \int_0^1 - \frac{t}{3} ( \int_{- t}^t F^{(4)}) - 60 t^2 G (1)
    \nonumber\\
    & \geqslant & \int_0^1 - \frac{t}{3} 2 tN - 60 t^2 G (1)
    \nonumber\\
    & = & - \frac{2}{3} (N + 90 G (1)) \int_0^1 t^2 \nonumber\\
    & = & - \frac{2}{3} (N + 90 G (1)) \frac{1}{3} . \nonumber
  \end{eqnarray}
  Hence, $N \geqslant - 90 G (1)$. Similarly, $0 \leqslant - \frac{2}{9} (- N
  + 90 G (1))$. Consequently, $90 G (1) \leqslant N$. We conclude that $|90G
  (1) | \leqslant N$.

A similar argument works to justify e.g.\ Romberg's integration
method~\cite{SuliMayers} which generalizes Simpson rule.
\section{Acknowledgements}

We thank Henri Lombardi for pointing us towards the Genocchi formula and we thank the referees for suggestions that
helped to improve the presentation of the article.

{The research leading to these results has received funding from the
European Union's 7th Framework Programme under grant agreement nr.\ 243847
(ForMath).}

\bibliographystyle{jloganal}

\end{document}